\documentclass[11pt,twoside,reqno]{amsart}
\usepackage[utf8]{inputenc}
\usepackage[T1]{fontenc}
\usepackage{import}
\usepackage{newclude}
\usepackage{amsmath,amssymb,amsthm,mathtools}

\usepackage{anysize}
\usepackage{amscd}
\usepackage{stmaryrd}
\usepackage{wasysym}
\usepackage{mathrsfs}
\usepackage[scr=boondox]{mathalfa}
\usepackage{enumitem}
\usepackage{accents}
\usepackage{dsfont}
\usepackage{soul}
\usepackage{color}
\usepackage[all]{xy}
\usepackage{float}
\usepackage{bm}
\usepackage{makecell}
\usepackage{thmtools}      
\usepackage{setspace}
\usepackage{comment}
\usepackage{tikz}
\usepackage{tikz-cd}
\usepackage{mathdots}
\usepackage{graphicx}
\usepackage[myheadings]{fullpage}

\usepackage[hypertexnames=false,colorlinks,linkcolor=black,citecolor=black]{hyperref}
\usepackage[capitalize,nameinlink]{cleveref}

\numberwithin{equation}{section}
\newcommand\mtop{.95in}
\newcommand\mbottom{.95in}
\newcommand\mleft{1in}
\newcommand\mright{1in}
\usepackage[top = \mtop, bottom = \mbottom, left = \mleft, right=\mright]{geometry}
\DeclareMathOperator{\val}{val}

\newtheorem{thm}{Theorem}[section]

\newtheorem{prop}[thm]{Proposition}

\newtheorem{lemma}[thm]{Lemma}

\theoremstyle{definition}
\newtheorem{defi}[thm]{Definition}

\newtheorem{rmk}[thm]{Remark}

\newcommand\reallywidehat[1]{%
\savestack{\tmpbox}{\stretchto{%
  \scaleto{%
    \scalerel*[\widthof{\ensuremath{#1}}]{\kern-.6pt\bigwedge\kern-.6pt}%
    {\rule[-\textheight/2]{1ex}{\textheight}}
  }{\textheight}%
}{0.5ex}}%
\stackon[1pt]{#1}{\tmpbox}%
}
\DeclareSymbolFont{bbold}{U}{bbold}{m}{n}
\DeclareSymbolFontAlphabet{\mathbbold}{bbold}

\makeatletter
\def\@tocline#1#2#3#4#5#6#7{\relax
  \ifnum #1>\c@tocdepth 
  \else
    \par \addpenalty\@secpenalty\addvspace{#2}%
    \begingroup \hyphenpenalty\@M
    \@ifempty{#4}{%
      \@tempdima\csname r@tocindent\number#1\endcsname\relax
    }{%
      \@tempdima#4\relax
    }%
    \parindent\z@ \leftskip#3\relax \advance\leftskip\@tempdima\relax
    \rightskip\@pnumwidth plus4em \parfillskip-\@pnumwidth
    #5\leavevmode\hskip-\@tempdima
      \ifcase #1
       \or\or \hskip 1em \or \hskip 2em \else \hskip 3em \fi%
      #6\nobreak\relax
    \hfill\hbox to\@pnumwidth{\@tocpagenum{#7}}\par
    \nobreak
    \endgroup
  \fi}
\makeatother

\newcommand{\R}{\mathbb{R}}
\newcommand{\Z}{\mathbb{Z}}
\newcommand{\Q}{\mathbb{Q}}
\newcommand{\N}{\mathbb{N}}

\newcommand{\C}{\mathbb{C}}

\newcommand{\D}{\mathbb{D}}
\newcommand{\I}{\mathcal{I}}
\newcommand{\U}{\mathbb{U}}

\newcommand{\bbone}{\mathbbold{1}}

\newcommand{\eps}{\epsilon}

\newcommand{\bad}{\text{bad}}

\renewcommand{\L}{\mathcal{L}}
\newcommand{\dist}{\text{dist}}

\newcommand{\T}{\mathcal{T}}

\renewcommand{\P}{\mathcal{P}}
\renewcommand{\D}{\mathcal{D}}

\title{High-low method and $p$-adic Furstenberg set over the plane}
\author{Kevin Ren and Jiahe Shen}

\date{\today}
\begin{document}

\thanks{We thank Ruixiang Zhang for encouraging us to invest in this project. We also thank Yuhan Chu for sharing her draft that introduces \Cref{prop: High-low method} to us. J.S. thanks his advisor Ivan Corwin for providing funding support with his NSF grant DMS-2246576 and Simons Investigator grant 929852. K.R. is supported by an NSF GRFP, Simons Dissertation Fellowship, and Cubist Fellowship.}

\maketitle

\begin{abstract} 
We establish a $p$-adic analogue of a recent significant result of Ren-Wang \cite{ren2023furstenberg} on Furstenberg sets in the Euclidean plane. Building on the $p$-adic version of the high-low method from Chu \cite{chu2025point}, we analyze cube-tube incidences in $\mathbb{Q}_p^2$ and prove that for $s < t < 2 - s$, any semi-well-spaced $(s,t)$-Furstenberg set over $\mathbb{Q}_p^2$ has Hausdorff dimension  $\ge\frac{3s+t}{2}$. Moreover, as a byproduct of our argument, we obtain the sharp lower bounds $s+t$ (for $0<t\le s\le 1$) and $s+1$ (for $s+t\ge 2$) for general $(s,t)$-Furstenberg sets without the semi-well-spaced assumption, thereby confirming that all three lower bounds match those in the Euclidean case.
\end{abstract}

\textbf{Keywords: }\keywords{$p$-adic Furstenberg sets; Hausdorff dimension; geometric measure theory}

\textbf{Mathematics Subject Classification (2020): }\subjclass{28A78 (primary); 28A80 (secondary)}

\tableofcontents

\section{Introduction}\label{sec: Introduction}

\subsection{A brief history over the Euclidean plane}

Fix $s\in(0,1]$ and $t\in(0,2]$. An $(s,t)$-Furstenberg set is a set $E\subset \R^2$ such that there exists a family $\L$ of lines with $\dim_H\L\ge t$ such that $\dim_H(E\cap l)\ge s$ for all $l\in\L$. The longstanding Furstenberg set conjecture asserts that when $E$ is a $(s,t)$-Furstenberg set, we have the sharp bound
$$\dim_HE\ge\min\{s+t,\frac{3s+t}{2},s+1\}.$$
The Furstenberg set conjecture can be regarded as the continuous analogue of the Szemerédi-Trotter theorem \cite[Theorem 1]{szemeredi1983extremal}. The bound $s+t$ for $0<t\le s\le 1$ was obtained by Lust-Stull \cite{lutz2020bounding} and independently by Héra-Shmerkin-Yavicoli \cite{hera2021improved}, while the bound $s+1$ for $s+t\ge 2$ was established by Fu-Ren \cite{fu2024incidence}. The product of a $t$-dim Cantor set and an $s$-dim Cantor set gives a sharp construction for the $s+t$ case, and the product of a line and an $s$-dim Cantor set gives a sharp construction for the $s+1$ case.

The intermediate regime $\frac{3s+t}{2}$ for $s<t<2-s$ is the most challenging case. Wolff \cite{wolff1999recent} constructed a sharp construction based on a grid-like configuration. Orponen-Shmerkin \cite{orponen2023projections} resolved the case in which the set is almost Ahlfors-David regular (that is, it exhibits comparable scaling behavior at all scales), using tools from projection theory and additive combinatorics. More recently, Ren-Wang \cite{ren2023furstenberg} completed the proof of the Furstenberg set conjecture in $\R^2$. Their work consists of two major components: first, they established the semi-well-spaced case via the high-low method in Fourier analysis, following ideas originating from Guth-Solomon-Wang \cite{guth2019incidence}; second, they used a multiscale decomposition technique \cite[Lemma 6.6]{ren2023furstenberg}, which partitions the scales into intervals such that on each interval the configuration is either almost AD-regular or semi-well-spaced. Since both scenarios are already resolved, the proof of the lower bound in the Euclidean plane follows.

\subsection{The $p$-adic case}

Let $p$ be a prime number. The field of $p$-adic numbers $\Q_p$ is the completion of $\Q$ with respect to the $p$-adic norm $||x||=p^{-\val_p(x)}$. This metric induces the usual notion of Hausdorff dimension for subsets $E\subset\Q_p^n$.

Motivated by the real case, it is natural to study point-line incidence in the $p$-adic setting. The Szemerédi-Trotter theorem in $\Q_p^2$ was proved in \cite[Theorem 1]{shen2025szemer}, and the expression is the same as the Euclidean plane. This makes us more willing to believe that the Hausdorff dimension of a $(s,t)$-Furstenberg set $E\subset\Q_p^2$ should have the sharp bound
$$\dim_H E\ge\min\{s+t,\frac{3s+t}{2},s+1\}.$$ 
Similar to the Euclidean strategy, a potential roadmap towards the proof could be the following:
\begin{enumerate}
\item Establishing the semi-well-spaced case via a $p$-adic version of high-low method using discrete Fourier analysis. \label{item: semi-well spaced step}
\item Proving the almost Ahlfors-David regular case using projection theory and additive combinatorics. \label{item: almost AD-regular step}
\item Partitioning the scales into intervals so that on each interval one reduces to either \eqref{item: semi-well spaced step} or \eqref{item: almost AD-regular step}. \label{item: partition scales step}
\end{enumerate}

Step \eqref{item: partition scales step} can be carried out along the lines of \cite[Lemma 6.6]{ren2023furstenberg} and is comparatively straightforward. Step \eqref{item: almost AD-regular step} requires asymmetric sum-product estimates in the ring $\Z/p^n\Z$, which are of independent interest. In this paper we focus on step \eqref{item: semi-well spaced step}, and the main results are as follows.

\begin{thm}\label{thm: lower bound of dimension for semi-well-spaced}
Let $E\subset\Q_p^2$ be an $(s,t)$-Furstenberg set with $s<t<2-s$, associated with a semi-well-spaced line set. Then we have
$$\dim_HE\ge\frac{3s+t}{2}.$$
\end{thm}

Here the term \emph{semi-well-spaced} refers to quantitative spacing assumptions. The reader may turn to \Cref{thm: discretized s<t<2-s} for a detailed description.

\begin{thm}\label{thm: sharp bound s+t}
Let $E\subset\Q_p^2$ be an $(s,t)$-Furstenberg set, where $0<t\le s\le 1$. Then we have the sharp bound $\dim_HE\ge s+t$.
\end{thm}

\begin{thm}\label{thm: sharp bound s+1}
Let $E\subset\Q_p^2$ be an $(s,t)$-Furstenberg set, where $s+t\ge 2$. Then we have the sharp bound $\dim_HE\ge s+1$.
\end{thm}

The sharp examples for \Cref{thm: sharp bound s+t} and \Cref{thm: sharp bound s+1} are analogous to the Euclidean constructions. For the $s+t$ case, take the product of a $t$-dim Cantor set and a $s$-dim Cantor set; for the $s+1$ case, take the product of a line and a $s$-dim Cantor set. Accordingly, the main task is to prove the stated inequalities, which will be carried out in \Cref{sec: Proof of main theorems}. As far as we know, prior to this work, there is no known improvement of the $p$-adic Furstenberg set problem beyond the classical bounds of Wolff \cite{wolff1999recent}, which holds for any local field. Our results thus demonstrate that the high-low method continues to be a powerful and effective tool outside the Euclidean setting. We hope that future work will address the almost AD-regular case, leading to a full resolution of the Furstenberg set problem over $\Q_p^2$.

\begin{rmk}
Research on Hausdorff and packing dimensions in fractal geometry has evolved significantly over the past few decades. Beyond the Furstenberg set problem considered in this paper, the major recent breakthrough was achieved by Wang-Zahl \cite{wang2025volume}, who building on their earlier works \cite{wang2022sticky,wang2025assouad}, proved that every Kakeya set in $\R^3$ has full Hausdorff and Minkowski dimension $3$. Their proof combines a multiscale induction scheme with tools from projection theory and additive combinatorics.

Meanwhile, Arsovski \cite{arsovski2021p} proved that every Kakeya set in $\Q_p^n$ has Hausdorff and Minkowski dimension $n$. Then by \cite[Remark 2]{arsovski2021p}, the result extends to Kakeya sets in $K^n$, where $K/\Q_p$ is an arbitrary finite extension. Unlike the Euclidean case, the proof of the $p$-adic case is relatively short and relies on a discrete valuation generalization of Dvir’s polynomial method from finite fields, see \cite{dvir2009size}. This suggests that, with respect to Kakeya-type problems, the $p$-adic field $\Q_p$ behaves more analogously to finite fields than to $\R$.

However, for $p$-adic Furstenberg sets, the situation is markedly different, as our approach is inspired more by the real field setting than by the finite field methods. As we will demonstrate, the multiscale induction and high-low method in \cite{ren2023furstenberg} for the Euclidean plane extends straightforwardly to the $p$-adic setting. We hope to see further exploration of the parallels between these different field settings in the future.
\end{rmk}

\subsection{Plan of the paper}

In \Cref{sec: Preliminaries}, we introduce $p$-adic cubes, tubes, and the relevant dimensional notions. We also quote a $p$-adic high-low method for estimating cube-tube incidences. In \Cref{sec: Pigeonholing}, we establish several pigeonholing and multiscale lemmas that will be involved in the proof of the main theorems. In \Cref{sec: Proof of main theorems}, we complete the argument by proving the results stated in \Cref{sec: Introduction}.

\section{Preliminaries}\label{sec: Preliminaries}

Throughout this paper, the letter $C$ will often denote a constant, while $K$ usually refers to a parameter depending on $\delta$.

We write $A\lesssim B$ if $A\le CB$ for some constant $C>0$, and the notation $A\gtrsim B$, $A\sim B$ follows in a similar way. We write $A\lesssim_n B$ to emphasize that the implicit constant may depend on $n$. The notation $\gtrapprox_\delta$ means $\gtrsim_N\log(1/\delta)^{C_N}$ for some number $N$ independent of $\delta$, and likewise $\lessapprox_\delta,\approx_\delta$.

The notation $|A|$ denotes either the cardinality of $A$ when $A$ is a finite set, or the Haar measure when $A$ is measurable and infinite. 

\begin{defi}
($p$-adic $\delta$-cubes) Let $x=(a,b)\in\Z_p^2$, $\delta=p^{-n}$ where $n\in\N$. Then the \emph{$\delta$-cube} centered at $x$ is defined by
$$B(x,\delta):=(a+p^n\Z_p,b+p^n\Z_p)=\{(a_0,b_0): ||a-a_0||\le\delta,||b-b_0||\le\delta\}.$$
\end{defi}

In contrast to the real setting, where $B(x,r)$ typically denotes a Euclidean ball, in the $p$-adic setting balls and cubes coincide due to the ultrametric structure. Hence we use the same notation $B(x,\delta)$. Moreover, any two $p$-adic cubes are either disjoint or one contains the other, reflecting the underlying tree-like structure of $\Z_p^2$, which often simplifies multiscale arguments.

For $n\in\Z_{\ge 0}$ and $\delta=p^{-n}$, denote by
$$\D_\delta(\Z_p^2):=\{B(x,\delta)\mid x\in\Z_p^2\}$$
as the set of $\delta$-cubes in $\Z_p^2$. Via the natural quotient map $\Z_p\rightarrow \Z_p/p^k\Z_p$, we can also identify $\D_\delta(\Z_p^2)$ with $(\Z/p^n\Z)^2$. For distinct cubes $\mathscr{p}_1=B((a_1,b_1),\delta),\mathscr{p}_2=B((a_2,b_2),\delta)\in\D_\delta$, we define their distance by
$$\dist(\mathscr{p}_1,\mathscr{p}_2):=\max\{||a_1-a_2||,||b_1-b_2||\},$$
which takes values in $\{p^{-k}:0\le k\le n-1\}$.

For a set $P\subset\Z_p^2$, denote
$$\D_\delta(P):=\{\mathscr{p}\in\D_\delta(\Z_p^2): \mathscr{p}\cap P\ne\emptyset\}.$$
For brevity, we write $\P_\delta=\D_\delta(\Z_p^2)$. If $\delta<\Delta\in p^{-\N}$, and $\mathscr{p}=B(x,\delta)\in\P_\delta$, let $\mathscr{p}^\Delta=B(x,\Delta)$ denote the unique $\Delta$-cube $\mathbf{p}\in\P_\Delta$ that contains $\mathscr{p}$.

\begin{defi}\label{defi: covering number}
For any set $P\subset\Z_p^2$ and $\delta=p^{-n}\in p^{-\N}$, the $\delta$-covering number of $P$ is defined by
$$|P|_\delta:=|\D_\delta(P)|.$$
\end{defi}

\begin{defi}\label{defi: dsC set}
($(\delta,s,C)$-set) For $\delta=p^{-n}\in p^{-\N},s\in[0,2]$ and $C>0$, a nonempty bounded set $A\subset\Z_p^2$ is called a \emph{$(\delta,s,C)$-set} if for all $r=p^{-m}\in[\delta,1]$ and $x\in\Z_p^2$,
$$|A\cap B(x,r)|_\delta\le Cr^s|A|_\delta.$$
\end{defi}

If $A$ is a $(\delta,s,C)$-set, then $|A|_\delta\ge C^{-1}\delta^{-s}$. If $\P\subset\P_\delta$ is a set of $\delta$-cubes, we say that $\P$ is a $(\delta,s,C)$-set if the union of $\delta$-cubes in $\P$ is a $(\delta,s,C)$-set in the above sense and we write $|\P|$ to denote the number of $\delta$-cubes in $\P$ and $|\P|_\Delta:=|\D_\Delta(\P)|$, so that $|\P|_\delta=|\P|$. 

\Cref{defi: dsC set} is to be contrasted with the following definition.

\begin{defi}
($(\delta,s,C)$-Katz-Tao-set) For $\delta=p^{-n}\in p^{-\N},s\in[0,2]$ and $C>0$, a nonempty bounded set $A\subset\Z_p^2$ is called a \emph{$(\delta,s,C)$-Katz-Tao set} if for all $r=p^{-m}\in[\delta,1]$ and $x\in \Z_p^2$,
$$|A\cap B(x,r)|_\delta\le C\left(\frac{r}{s}\right)^\delta.$$
\end{defi}

\begin{defi}
($p$-adic $\delta$-tubes) Let $\delta=p^{-n}\in p^{-\N}$. A $\delta$-tube is a set of the form $T=\cup_{x\in\mathscr{p}}\mathbf{D}(x)$, where $\mathscr{p}\in\P_\delta$, and $\mathbf{D}$ is the point-line duality map
$$\mathbf{D}(a,b)=\{(x,y)\in\Z_p^2: y=ax+b\}\subset \Z_p^2$$
sending the point $(a,b)$ to the $p$-adic line with slope $a$ and intercept $b$. For brevity, we write $\mathbf{D}(\mathscr{p}):=\cup_{x\in\mathscr{p}}\mathbf{D}(x)$ as the $\delta$-tube that corresponds to $\mathscr{p}$. The collection of all $p$-adic $\delta$-tubes is denoted
$$\T_\delta:=\{\mathbf{D}(\mathscr{p}): \mathscr{p}\in\P_\delta\},$$
which may also be viewed as the set of lines in $(\Z/p^n\Z)^2$ with slopes in $\Z/p^n\Z$.
\end{defi}

A collection of $p$-adic $\delta$-tubes $\{\mathbf{D}(\mathscr{p})\}_{\mathscr{p}\in\P}$ is called a $(\delta,s,C)$-set if $\P$ is a $(\delta,s,C)$-set. For two distinct cubes $\mathscr{p}_1,\mathscr{p}_2\in\P_\delta$ with $$\dist(\mathscr{p}_1,\mathscr{p}_2)=p^{-k},\quad 0\le k\le n-1,$$ 
there are at most $p^k$ lines in $(\Z/p^n\Z)^2$ that pass through both of them.

If $\delta<\Delta\in p^{-\N}$ and $T\in\T_\delta$, denote by $T^\Delta$ the unique $p$-adic tube $\mathbf{T}\in\T_\delta$ containing $T$. For a set of $p$-adic $\delta$-tubes $\T$ and any scale $\Delta>\delta$, define the $p$-adic $\Delta$-covering number $|\T|_\Delta:=|\D_\Delta(\T)|$, where $\D_\Delta(\T):=\{\mathbf{T}\in\T_\Delta:\exists T\in\T,\mathbf{T}=T^\Delta\}$. We also denote $\T\cap\mathbf{T}:=\{T\in\T: T\subset\mathbf{T}\}$ for all $\mathbf{T}\in\D_\Delta(\T)$. When $\delta=\Delta$, we usually omit the subscript and simply write $|\T|=|\T|_\delta$.


\subsection{Fourier transform}

In this subsection, we fix $n\in\N$. For a function $f: (\Z/p^n\Z)^2\rightarrow\C$, its Fourier transform $\widehat f:(\Z/p^n\Z)^2\rightarrow\C$ is defined by
$$\widehat f(\xi):=p^{-n}\sum_{x\in(\Z/p^n\Z)^2}f(x)\exp(-\frac{2\pi i\langle\xi,x\rangle}{p^n}).$$
Here $\langle\cdot,\cdot\rangle$ denotes the standard dot product on $(\Z/p^n\Z)^2$, given by 
$$\langle (a_1,b_1),(a_2,b_2)\rangle=a_1a_2+b_1b_2\mod p^n.$$

With this normalization, the Fourier transform is an involution, i.e., $\hat{\hat f}=f$. The following propositions are $p$-adic analogues of standard results in Euclidean Fourier analysis; their proofs are routine and therefore omitted.

\begin{prop}
(Plancherel/Parseval identity) For any functions $f,g: (\Z/p^n\Z)^2\rightarrow\C$, We have the equality
\begin{equation}\label{eq: Plancherel}
\sum_{x\in(\Z/p^n\Z)^2}f(x)\bar g(x)=\sum_{\xi\in(\Z/p^n\Z)^2}\widehat f(\xi)\bar{\widehat g}(\xi).
\end{equation}
In particular, when $f=g$, we have
\begin{equation}\label{eq: Plancherel_f=g}
\sum_{x\in(\Z/p^n\Z)^2}|f(x)|^2=\sum_{\xi\in(\Z/p^n\Z)^2}|\widehat f(\xi)|^2.
\end{equation}
\end{prop}

\begin{prop}
(Convolution identity) For any functions $f,g: (\Z/p^n\Z)^2\rightarrow\C$, We have the equality
\begin{equation}
\widehat{f*g}(\xi)=p^n\widehat f(\xi)\widehat g(\xi).
\end{equation}
Here the convolution function $f*g:(\Z/p^n\Z)^2\rightarrow\C$ is defined as
$$(f*g)(x)=\sum_{y\in(\Z/p^n\Z)^2}f(y)g(x-y).$$
\end{prop}

\subsection{High-low method for incidences} 

\begin{defi}
(Incidence count) Let $\delta<\Delta\in p^{-\N}$. Let $\P$ be a finite collection of distinct $\delta$-cubes equipped with a weight function $w:\P\rightarrow\R_{\ge 0}$ and let $\T$ be a set of (not necessarily distinct) $\Delta$-tubes. We define the weighted incidence number between $\P$ and $\T$ by
$$\I_w(\P,\T):=\sum_{\mathscr{p}\in\P}\sum_{T\in\T,\mathscr{p}\subset T}w(\mathscr{p}).$$
If $w=1$, we often omit the subscript and write $\I(\P,\T):=\I_w(\P,\T)$.
\end{defi}

The following proposition is a slightly different version of \cite[Proposition 3.4]{chu2025point}, which is a $p$-adic analogue of the high-low method in \cite[Section 2]{guth2019incidence}.

\begin{prop}\label{prop: High-low method}
Fix $\delta=p^{-n}\in p^{-\N}$. Let $\P$ be a set of distinct $p$-adic $\delta$-cubes contained in $\Z_p^2$ with a weighted function $w: \P\rightarrow\R^{\ge 0}$ and $\T$ a set of distinct $p$-adic $\delta$-tubes in $\Z_p^2$. Then for all $1\le k\le n-1$, we have
$$\I_w(\P,\T)\le p^{\frac{n+k-1}{2}}|\T|^{1/2}\left(\sum_{\mathscr{p}\in\P}w(\mathscr{p})^2\right )^{1/2}+p^{-k}I_w(\P,\T^{p^k\delta}).$$
Here $\T^{p^k\delta}:=\{T^{p^k\delta}:T\in\T\}$ is counted with multiplicity.
\end{prop}

\begin{proof}[Proof of \Cref{prop: High-low method}]
Define $f,g: (\Z/p^n\Z)^2\rightarrow\C$ by
$$f(x)=\sum_{\mathscr{p}\in\P}w(\mathscr{p})\bbone_{\mathscr{p}}(x),\quad g(x)=\sum_{T\in\T}\bbone_{T}(x),$$
where 
$$\bbone_{\mathscr{p}}(x)=\begin{cases}
1 & \mathscr{p}=B(x,\delta) \\
0 & \text{otherwise}
\end{cases},\quad\bbone_{T}(x)=\begin{cases}
1 & x\text{ lies in }T\\
0 & \text{otherwise}
\end{cases}.$$
Then by definition, we have
$$\I_w(\P,\T)=\sum_{x\in(\Z/p^n\Z)^2}f(x)g(x)=\sum_{\xi\in(\Z/p^n\Z)^2}\widehat f(\xi)\bar{\widehat g}(\xi).$$
Let $\eta=\bbone_{B(0,p^{-k})}$. We decompose the above sum as
\begin{equation}\label{eq: high-low decomposition}
\I_w(\P,\T)=\sum_{\xi\in(\Z/p^n\Z)^2}\widehat f(\xi)\bar{\widehat g}(\xi)\eta(\xi)+\widehat f(\xi)\bar{\widehat g}(\xi)(1-\eta(\xi)):=L+H,
\end{equation}
where 
$$L:=\sum_{\xi\in(\Z/p^n\Z)^2}\widehat f(\xi)\bar{\widehat g}(\xi)\eta(\xi)$$ 
is the \emph{low-frequency term}, and $$H:=\sum_{\xi\in(\Z/p^n\Z)^2}\widehat f(\xi)\bar{\widehat g}(\xi)(1-\eta(\xi))$$ 
is the \emph{high-frequency term}. Now let us estimate the two terms separately. For the high-frequency term, applying the Cauchy-Schwarz inequality,
\begin{align}\label{eq: Cauchy-Schwartz over f and g}
\begin{split}
H&=\sum_{\xi\in(\Z/p^n\Z)^2}\widehat f(\xi)\bar{\widehat g}(\xi)(1-\eta(\xi))\\
&\le \left(\sum_{\xi\in(\Z/p^n\Z)^2}|\widehat f(\xi)|\right)^{1/2}\cdot\left(\sum_{\xi\in(\Z/p^n\Z)^2}|{\widehat g}(\xi)|^2(1-\eta(\xi))^2\right)^{1/2}\\
&=\left(\sum_{x\in(\Z/p^n\Z)^2}| f(x)|\right)^{1/2}\cdot\left(\sum_{\xi\in(\Z/p^n\Z)^2}|{\widehat g}(\xi)|^2(1-\eta(\xi))^2\right)^{1/2}\\
&=\left(\sum_{\mathscr{p}\in\P}| w(\mathscr{p})|^2\right)^{1/2}\cdot\left(\sum_{\xi\in(\Z/p^n\Z)^2}|{\widehat g}(\xi)|^2(1-\eta(\xi))^2\right)^{1/2}.
\end{split}
\end{align}
Here the third line comes from Parseval's identity \eqref{eq: Plancherel_f=g}.

Observe that for all $T\in\T$, $\widehat\bbone_{T}$ is supported on the line $l'$ through $(0,0)$ and perpendicular to $T$. For all $\xi\in(\Z/p^n\Z)^2\backslash(p^k\Z/p^n\Z)^2$, there are at most $p^{k-1}$ distinct lines in $(\Z/p^n\Z)^2$ passing through $0$ and $\xi$. Denote by $\mathbb{P}(\Z/p^n\Z)^2$ the set of directions of lines in $(\Z/p^n\Z)^2$. Applying Cauchy-Schwarz, we have 
\begin{align}\label{eq: estimate of g}
\begin{split}
\sum_{\xi\in(\Z/p^n\Z)^2}|\widehat g(\xi)|^2(1-\eta(\xi))^2&=\sum_{\xi\in(\Z/p^n\Z)^2}(1-\eta(\xi))^2|\sum_{T\in\T}\widehat\bbone_T(\xi)|^2\\
&\le p^{k-1}\sum_{\xi\in(\Z/p^n\Z)^2}(1-\eta(\xi))^2\sum_{d\in\mathbb{P}(\Z/p^n\Z)^2}|\sum_{T\text{ has direction }d}\widehat\bbone_T(\xi)|^2\\
&\le p^{k-1}\sum_{\xi\in(\Z/p^n\Z)^2,d\in\mathbb{P}(\Z/p^n\Z)^2}|\sum_{T\text{ has direction }d}\widehat\bbone_T(\xi)|^2.\\
\end{split}
\end{align}
Applying \eqref{eq: Plancherel_f=g}, we have
\begin{align}\label{eq: decomposition of direction and Parseval}
\begin{split}
\sum_{\xi\in(\Z/p^n\Z)^2,d\in\mathbb{P}(\Z/p^n\Z)^2}|\sum_{T\text{ has direction }d}\widehat\bbone_T(\xi)|^2&=\sum_{d\in\mathbb{P}(\Z/p^n\Z)^2}\sum_{x\in(\Z/p^n\Z)^2}|\sum_{T\text{ has direction }d}\bbone_T(x)|^2\\
&=\sum_{d\in\mathbb{P}(\Z/p^n\Z)^2}p^k|T\in\T: T\text{ has direction }d|\\
&=p^{n+k-1}|\T|.
\end{split}
\end{align}
Here the second line holds because parallel lines are disjoint. The result \eqref{eq: Cauchy-Schwartz over f and g}, \eqref{eq: estimate of g} and \eqref{eq: decomposition of direction and Parseval} together yields
\begin{equation}\label{eq: estimate of high}
H\le p^{\frac{n+k-1}{2}}|\T|^{1/2}\left(\sum_{\mathscr{p}\in\P}w(\mathscr{p})^2\right )^{1/2}.
\end{equation}
For the low-frequency term, applying \eqref{eq: Plancherel} and we have
\begin{align}\label{eq: L as convolution}
\begin{split}
L&=\sum_{\xi\in(\Z/p^n\Z)^2}\widehat f(\xi)\bar{\widehat g}(\xi)\eta(\xi)\\
&=p^{-n}\sum_{x\in(\Z/p^n\Z)^2} f(x)(g*h)(x).\\
\end{split}
\end{align}
Here $h=p^{n-2k}\bbone_{B(0,p^{k-n})}=\widehat\eta$ is the Fourier transform of $\eta$. We now study the convolution $g*h$,
\begin{align}\label{eq: convolution of g and h}
\begin{split}
(g*h)(x)&=\sum_{y\in(\Z/p^n\Z)^2}g(y)h(x-y)\\
&=p^{n-2k}\sum_{y\in(\Z/p^n\Z)^2}\sum_{T\in\T}\bbone_{T}(y)\bbone_{B(x,p^{k-n})}(y)\\
&=p^{n-2k}\sum_{T\in\T}|T\cap B(x,p^{k-n})|\\
&=p^{n-k}|\{T\in\T: T\cap B(x,p^{k-n})\ne \emptyset\}|\\
&=p^{n-k}|\{T\in\T: T^{p^k\delta}\supset B(x,p^{-n})\}|.
\end{split}
\end{align}
The third line holds because for each $T\in\T$ such that $|T\cap B(x,p^{k-n})|>0$, we must have $|T\cap B(x,p^{k-n})|=p^k$. Combining the result in \eqref{eq: L as convolution} and \eqref{eq: convolution of g and h}, we have
\begin{align}\label{eq: equality of low}
\begin{split}
L&=p^{-k}\sum_{x\in(\Z/p^n\Z)^2}f(x)\cdot|\{T\in\T: T^{p^k\delta}\supset B(x,p^{-n})\}|\\
&=p^{-k}\sum_{\mathscr{p}\in\P}\sum_{T\in\T,\mathscr{p}\subset T^{p^k\delta}}w(\mathscr{p})\\
&=p^{-k}\mathcal{I}_w(\P,\T^{p^k\delta}).
\end{split}
\end{align}
The decomposition \eqref{eq: high-low decomposition}, together with the results \eqref{eq: estimate of high} and \eqref{eq: equality of low} gives the proof.
\end{proof}

\section{Pigeonholing}\label{sec: Pigeonholing}

In this section, we make an analogy of \cite[Section 3]{ren2023furstenberg} and introduce some pigeonholing that enables us to reduce our study to circumstances with good properties. 

\begin{defi}
Fix $\delta=p^{-n}\in p^{-\N},s\in[0,1],C>0$, and $M\in\N$. We say that a pair $(\P,\T)\in\P_\delta \times \T_\delta$ is a \emph{$\delta$-configuration} if every $\mathscr{p}\in\P$ is associated with a subset $\T(\mathscr{p})\subset\T$ such that $\mathscr{p}\subset T$ for all $T\in\T(\mathscr{p})$ (note that $\T(\mathscr{p})$ does not need to be all the tubes in $\T$ that contain $\mathscr{p}$). We say that the pair $(\P,\T)$ is a \emph{$(\delta,s,C,M)$-nice configuration} if for every $\mathscr{p}\in\P$, $\T(\mathscr{p})\subset\T$ is a $(\delta,s,C)$-set with $|\T(p)|\sim M$.
\end{defi}

It follows from the basic property of the $(\delta,s,C)$-set that in a $(\delta,s,C,M)$-nice configuration, we have $M\gtrsim C^{-1}\delta^{-s}$.

\begin{defi}
Let $\P\subset\P_\delta$. We say that the subset $\P'\subset \P$ is a \emph{refinement} of $\P$ if $|\P'|_\delta\approx_\delta|\P|_\delta$. Let $\Delta=p^{-m}\in(\delta,1)$, we say that $\mathcal{P}'\subset \P$ is a \emph{refinement of $\P$ at resolution $\Delta$} if there is a refinement $\P_\Delta'$ of $\D_\Delta(\P)$ such that $\P'=\bigcup_{\mathbf{p}\in\P_\Delta'}(\P\cap\mathbf{p})$ (note that in this case, $\P'$ does not necessarily satisfy $|\P'|_\delta\approx_\delta|\P|_\delta$.). We define the refinement over $p$-adic $\delta$-tubes $\T$ at resolution $\Delta$ similarly.
\end{defi}

\begin{defi}
Let $(\P_0,\T_0)$ be a $(\delta,s,C_0,M_0)$-configuration. We say that a $\delta$-configuration $(\P,\T)$ is a refinement of $(\P_0,\T_0)$ if $\P\subset\P_0$ if $\P\subset \P_0$, $|\P|_\delta\approx_\delta|\P_0|_\delta$ (i.e, $\P$ is a refinement at $\P_0$), and for any $\mathscr{p}\in\P$, there is $\T(\mathscr{p})\subset\T_0(\P)\cap\T$ with $\sum_{\mathscr{p}\in\P}|\T(\mathscr{p})|_\delta\gtrapprox_\delta |\P_0|_\delta\cdot M_0$. Furthermore, we say that a refinement $(\P,\T)$ of $(\P_0,\T_0)$ is a nice configuration refinement if $|\T(\mathscr{p})|\approx_\delta|\T_0(\mathscr{p})|$ for all $\mathscr{p}\in\P$.

Let $\Delta=p^{-m}\in(\delta,1)$, we say that $(\P,\T)$ is a refinement of $(\P_0,\T_0)$ at resolution $\Delta$ if $\P$ is a refinement of $\P_0$ at resolution $\Delta$ and for each $p\in\P$, $\mathbf{T}\in\D_\Delta(\T(p)),\mathbf{T}\cap\T(\mathscr{p})=\T\cap\T_0(\mathscr{p})$. 
\end{defi}

The following proposition claims that when making a refinement, we can always pigeonhole so that the refinement is a nice configuration refinement.

\begin{prop}
For any $(\delta,s,C_0,M_0)$-nice configuration $(\P_0,\T_0)$ and refinement $(\P,\T)$ of $(\P_0,\T_0)$, there exists a refinement $(\P',\T')$ of $(\P,\T)$ that is a nice configuration refinement of $(\P_0,\T_0)$.
\end{prop}

\begin{defi}
Let $N\ge 1$, and
$$\delta=\Delta_N<\Delta_{N-1}<\cdots<\Delta_1<\Delta_0=1$$
be a sequence of $p$-adic scales. We say that a set $P\subset\Z_p^2$ is $\{\Delta_j\}_{j=1}^N$-uniform if there is a sequence $\{N_j\}_{j=1}^N$ with $N_j\in p^{\N}$ and $|P\cap Q|_{\Delta_j}\in[N_j/p,N_j)$ for all $j=\{1,2,\ldots,N\}$ and $Q\in\D_{\Delta_{j-1}}(\P)$.
\end{defi}

The following lemma is an analog of \cite[Lemma 3.6]{shmerkin2019furstenberg}, which states that given a set in $\Z_p^2$, we can always pigeonhole a relatively large subset that is uniform. The proof technique is the same, so we omit it here.

\begin{lemma}\label{lem: pigeonhole}
Let $P\in\Z_p^2$, $N,T\in\N$, and $\delta=p^{-NT}$. Also, let $\Delta_j:=p^{-jT}$ for all $0\le j\le N$. Then there exists a $\{\Delta_j\}_{j=1}^N$-uniform set $P'\subset P$ such that
$$|P'|_\delta\ge (pT)^{-N}|P|_\delta=p^{-\frac{\log_p T+1}{T}\cdot NT}|P|_\delta.$$
In particular, if $\epsilon>0$ and $T^{-1}(\log_p T+1)\le\epsilon$, then $|P'|_\delta\ge\delta^{\epsilon}|P|_\delta$.
\end{lemma}



\section{Proof of main theorems}\label{sec: Proof of main theorems}

The main goal of this section is to prove \Cref{thm: lower bound of dimension for semi-well-spaced}, \Cref{thm: sharp bound s+t} and \Cref{thm: sharp bound s+1} in \Cref{sec: Introduction}. In order to prove these theorems, we transfer them to the following dual discretized version. (This is a usual step for studying fractal sets, see \cite[Theorem 4.1]{ren2023furstenberg} over the real field case for instance.)

\begin{thm}\label{thm: discretized s<t<2-s}
Suppose $s\in(0,1)$, and $s<t<2-s$. Then for every $\epsilon>0$ and $\eta=\frac{\epsilon^2}{1225}$, the following holds for any $\delta\in p^{-\N}$ sufficiently small. Let $(\P,\T)$ be a $(\delta,s,\delta^{-\eta},M)$-nice configuration, so that $M\ge \delta^{-s+\eta}$. Suppose $|\P|_\delta\sim\delta^{-t}$, and $\P\subset\P_\delta$ satisfies the stronger spacing condition
\begin{equation}\label{eq: strong space condition}
|\P\cap Q|_\delta\lesssim\delta^{-\eta}\cdot\max(\rho^{2-s}|\P|_\delta,(\rho/\delta)^s),\quad\forall Q\in\D_\rho(\P),\delta\le\rho\le 1.
\end{equation}
Then we have $|\T|_\delta\gtrsim_\epsilon \delta^{-\frac{s+t}{2}+\epsilon}M.$
\end{thm}

\begin{thm}\label{thm: discretized 0<t<s<1}
Suppose $0<t\le s\le 1$. For every $\epsilon>0$, there exists $\eta>0$ such that the following holds for any $\delta\in p^{-\N}$ sufficiently small. Let $(\P,\T)$ be a $(\delta,s,\delta^{-\eta},M)$-nice configuration with $s\in(0,1]$, $\P$ be a $(\delta,t,\delta^{-\eta})$-set.
Then $|\T|_\delta\gtrsim_\epsilon\delta^{-t+\epsilon}M$.
\end{thm}

\begin{thm}\label{thm: discretized s+t>2}
Suppose $s+t\ge 2$. For every $\epsilon>0$, there exists $\eta>0$ such that the following holds for any $\delta\in p^{-\N}$ sufficiently small. Let $(\P,\T)$ be a $(\delta,s,\delta^{-\eta},M)$-nice configuration with $s\in(0,1]$, $\P$ be a $(\delta,t,\delta^{-\eta})$-set.
Then $|\T|_\delta\gtrsim_\epsilon\delta^{-1+\epsilon}M$.
\end{thm}

The proofs for \Cref{thm: discretized s<t<2-s}, \Cref{thm: discretized 0<t<s<1}, and \Cref{thm: discretized s+t>2} are all based on the following theorem, which is the analog of \cite[Proposition 4.6]{ren2023furstenberg}. It studies the key special case that $\P$ is $(2-s)$-dimensional at high scales and $s$-dimensional at low scales.

\begin{thm}\label{thm: 2-s high and s low}
Let $\delta,\Delta\in p^{-\N},\delta\le\Delta,\epsilon>0,s\in(0,1]$, and $(\P,\T)$ be a $(\delta,s,\delta^{-\epsilon^2},M)$-nice configuration, so that $M\gtrsim\delta^{-s+\epsilon^2}$. Suppose that $\P\in\D_\delta$ satisfies the spacing conditions
\begin{equation}\label{eq: 2-s high}
|\P\cap Q|_\delta\le\delta^{-\epsilon^2}\cdot\rho^{2-s}\cdot|\P|_\delta,\quad\forall\rho\in(\Delta,1),Q\in\D_\rho(\P),
\end{equation}
\begin{equation}\label{eq: s low}
|\P\cap Q|_\delta\le\delta^{-\epsilon^2}\cdot(\rho/\delta)^{s},\quad \forall\rho\in(\delta,\Delta),Q\in\D_\rho(\P).
\end{equation}
Then $|\T|_\delta\gtrapprox_{\delta,\epsilon}\delta^{35\epsilon}\cdot\min\{M|\P|_\delta,M^{3/2}|\P|_\delta^{1/2},\delta^{-1}M\}.$ Here the notation $\gtrapprox_{\delta,\epsilon}$ means $\le C_\epsilon\log(1/\delta)^C$, which is the same as the definition in \cite[Proposition 4.4]{ren2023furstenberg}.
\end{thm}

We will focus on the proof of \Cref{thm: 2-s high and s low} for the upcoming subsections and finally revisit \Cref{thm: discretized s<t<2-s}, \Cref{thm: discretized 0<t<s<1}, and \Cref{thm: discretized s+t>2} in \Cref{subsec: Returning to the discretized version of Furstenberg set} to complete the proof.

\subsection{Pigeonholing for building multiscales}

The discretized version of Furstenberg set conjecture enables us to start from a small fixed scale $\delta$, but it does not provide information for structures over the medium scales. Fortunately, we are able to build multiscale structures based on the pigeonholing results in \Cref{sec: Pigeonholing}. We will start from the following definition.

\begin{defi}
Fix $\delta\le\Delta\in p^{-\N}$, $s\in[0,1]$, and $C_\delta>0,M_\delta\in\N$. We say that a $\Delta$-configuration $(\mathcal{Q},\T^\Delta)$ covers a $(\delta,s,C_\delta,M_\delta)$-nice configuration $(\P,\T)$ if
$$\sum_{\mathscr{q}\in\mathcal{Q}}\sum_{\mathscr{p}\in\P}\sum_{\mathbf{T}\in\T^\Delta(\mathbf{p})}|\T_\mathscr{p}\cap\mathbf{T}|_\delta\gtrapprox_\delta|\P|_\delta\cdot M_\delta.$$
\end{defi}

The following proposition refines a nice configuration at scale $\delta$, which builds a medium scale $\Delta$ that also has nice configuration.

\begin{prop}\label{prop: new scale Delta}
Fix $\delta\in p^{-\N}$, $s\in(0,1]$, $M\in\N$ and let $\Delta\in(\delta,1)$. Let $(\P,\T)$ be a $(\delta,s,C,M)$-configuration, then there exists a refinement $(\P_1,\T_1)$ of $(\P,\T)$ such that $\P_1$ is a $\{1,\Delta,\delta\}$-uniform set and
\begin{enumerate}
\item $(\P_1^\Delta,\T_1^\Delta),\P_1:=\D_\Delta(\P_1)$ and for any $\mathbf{p}\in\P_1^\Delta,\T_1(\mathbf{p})^\Delta:=\D_\Delta(\cup_{p\in\P_1\cap\mathbf{p}}\T_1(\mathbf{p}))$, is a $(\Delta,s,C_\Delta,M_\Delta)$-nice configuration with $C_\Delta\lessapprox_\delta C$ and $M_\Delta\in\N$.
\item Any refinement $(\mathcal{Q},\T^\Delta)$ of $(\P_1^\Delta,\T_1^\Delta)$ covers $(\P,\T)$.
\end{enumerate}
\end{prop}

\begin{proof}
For every $\mathscr{p}\in \P$, applying \Cref{lem: pigeonhole} to the set $\T(\mathscr{p})$, we deduce that there exists a $\{1, \Delta, \delta\}$-uniform subset $\T_1(\mathscr{p})\subset \T(\mathscr{p})$ with $|\T_1(\mathscr{p})|\approx_{\delta} M$. Let $m(p)=|\T_1(\mathscr{p})|_{\Delta}$. Then for every $\mathbf{T}\subset\D_\Delta(\T)$ such that $\mathbf{T}\cap  \T_1(\mathscr{p}) \neq \emptyset$, we have $|\mathbf{T}\cap \T_1(\mathscr{p})|_{\delta} \approx_{\delta} \frac{M}{m(p)}$.
Applying \Cref{lem: pigeonhole} to the $\P$ again, there exists $m\in p^{\N}$ and a $\{1, \Delta, \delta\}$-uniform subset $\P_1\subset \P$, $|\P_1|_{\delta} \approx_{\delta} |\P|_{\delta}$    and  for every $\mathscr{p}\in \P_1$, $m\sim m(\mathscr{p})$.     

For each $\mathbf{p} \in \D_{\Delta}(\P_1)$, there exists $X(\mathbf{p})\in p^{\N}$ and a subset $\T_1(\mathbf{p})\subset \T_{\Delta}$ satisfying for each $\T\in \T_1(\mathbf{p})$, 
$$
|\{ p\in \P_1\cap \mathbf{p}: \T\in \D_{\Delta}(\T_1(p))\}|\sim X(\mathbf{p})
$$
and 
\begin{equation}\label{eq: X}
X(\mathbf{p}) \cdot |\T_1(\mathbf{p})|_{\Delta} \approx_{\Delta}|\P_1\cap \mathbf{p}|_{\delta} \cdot  m.
\end{equation}

By pigeonholing over $\D_\Delta(\P_1)$, there exists a common $X\in p^{\N}$ and a refinement of $\P_1$ at resolution $\Delta$, which we still denote $\P_1$, such that for each $\mathbf{p}\in \D_{\Delta}(\P_1)$, $X(\mathbf{p})=X$. For every $\mathscr{p}\in\P_1$, we replace $\T_1(\mathscr{p})$ by $\{T\in\T_1(\mathscr{p}): T^\Delta\in\T_1(\mathscr{p}^\Delta)\}$. Following the same proof as in \cite[Proposition 4.3]{ren2023furstenberg}, the refinement $(\P_1,\T_1)$ satisfies our requirements.
\end{proof}

Now we can iterate the construction in \Cref{prop: new scale Delta} to give multiscale structures to the configuration, which is what we do in the following proposition.

\begin{prop}\label{prop: multi-grid}
Fix $\Delta<\Delta_0<1,n\in\N,\epsilon\in(0,1/n]$, and let $\delta=\Delta^n$. Let $(\P_0,\T_0)$ be a $(\delta,s,C_0,M_0)$-nice configuration with $\P_0$ a $(\delta,t,C_0)$-set, $s+t\ge 2$. Then there exists a refinement $(\P,\T)$ of $(\P_0,\T_0)$ such that for any $w=\Delta^k,1\le k\le n$ and any $T\in\T$, 
\begin{equation}|T^w\cap\D_w(\P)|\lessapprox_{\delta,\epsilon} C_0^2\Delta^{-1}|\P|_w\cdot w^{1-\epsilon}.\label{eqn:tube est}\end{equation}
Here $|T^w\cap\D_w(\P)|$ means the number of $\mathbf{p}\in\D_w(\P)$ contained in $T^w$.
\end{prop}

\begin{proof}
The proof follows the same outline as \cite[Proposition 4.4]{ren2023furstenberg}, using Proposition \ref{prop: High-low method} as the key input for high-low method. So here we will sketch the key steps and then refer to \cite{ren2023furstenberg} for more exact details.

We perform an induction on $n$, the base case $n = 1$ being the geometric fact $|T \cap \D_\Delta (\P)| \le \Delta^{-1}$.

For the inductive step, assume true for $n-1$, and we will now prove the proposition statement for $n$. Without loss of generality, assume $\P$ is uniform at scales $\{ \Delta^i \}_{i=0}^n$. Apply Proposition \ref{prop: new scale Delta} with $\Delta^{n-1}$ in place of $\Delta$ to obtain a refinement $(\P_1, \T_1)$ of $(\P, \T)$. Now by the inductive hypothesis on $(\P_1^{\Delta^{n-1}}, \T_1^{\Delta^{n-1}})$ and a further refinement procedure found in \cite{ren2023furstenberg}, we can find a further refinement of $(\P_1, \T_1)$, which we call $(\P_2, \T_2)$, such that:
\begin{itemize}
    \item $|T \cap \P_2| \sim r$ where $r \in [1, \delta^{-1}]$;

    \item $|T^w\cap\D_w(\P)|\lessapprox_\delta C_0^2\Delta^{-1}|\P|_w\cdot w^{1-\epsilon}$ for any $w=\Delta^k,1\le k\le n-1$ and $T \in \T$.
\end{itemize}
Finally, we use Proposition \ref{prop: High-low method} to prove an upper bound on $r$ which matches \eqref{eqn:tube est} for $k = n$ (see \cite{ren2023furstenberg} for details).
\end{proof}

\subsection{Proof of \Cref{thm: 2-s high and s low}}

For any $T\in\T_\delta$ and $b\in\N$, let $N_{w,b}(T)$ be the number of $Q\in\D_w(\P)$ (count with multiplicity) that $|T\cap Q\cap \P|\ge b$. The following lemma is crucial for proving \Cref{thm: 2-s high and s low}.

\begin{lemma}\label{lem: large incidence tube}
Let $\delta\le\Delta\in p^{-\N},\epsilon\in\frac{1}{\N}$, and let $\P$ be a multi-set of $\delta$-cubes such that for all $1\le k\le \epsilon^{-1}$, each $Q\in\D_{\Delta^{k\epsilon}}(\P)$ contains about the same number of cubes in $\P$ (with multiplicity). For $a\ge 2$ and $ab>\delta^{1-2\epsilon}|\P|$, let $\T_{a,b}$ be a set of distinct $\delta$-tubes satisfying
\begin{enumerate}
\item $|T^w\cap\D_w(\P)|\le\Delta^{-\epsilon}\cdot w|\D_w(\P)|$ for all $w\in\{\Delta^\epsilon,\Delta^{2\epsilon},\ldots,\Delta\}$;
\item $N_{\Delta,b}(T)\ge a$.
\end{enumerate}
Then $|\T_{a,b}|\lesssim_\epsilon\frac{|\P|^2}{a^3b^2}\delta^{-5\epsilon}$.
\end{lemma}

We now explain how \Cref{lem: large incidence tube} implies \Cref{thm: 2-s high and s low}. The following lemma turns out to be a key pigeonholing.

\begin{lemma}\label{lem: many good tubes}
Fix $\delta<\Delta\in p^{-\N}$. Let $\P$ be a set of distinct $p$-adic $\delta$-tubes in $B(0,\Delta)$ such that $|\P\cap B(x,w)|\le K_1\cdot(\frac{w}{\delta})^s$ for all $\delta\le w\le \Delta$ and $w$-balls $B(x,w)$. For each $\mathscr{p}\in\P$, let $\T(\mathscr{p})$ be a $(\delta,s,K_2)$-set of $\delta$-tubes such that $\T(\mathscr{p})\sim M$. Let $\T=\cup_{\mathscr{p}\in\P}\T(\mathscr{p})$. Then there exists a subset $\P'\subset\P$ with $|\P'|\ge\frac{1}{2}|\P|$ and a subset $\T'(p) \subset \T(p)$ with $|\T'(p)| \ge \frac{1}{2} |\T(p)|$ such that for each $T\subset\T':=\cup_{\mathscr{p}\in\P'}\T'(\mathscr{p})$ we have $|\mathscr{p}\in\P':T\in\T'(\mathscr{p})|\le C_1K_1 K_2\log(1/\delta)$, for some universal constant $C_1>0$. 
\end{lemma}

\begin{proof}
Denote by $\T_{\bad}\subset\T$ the set of tubes with $|\{\mathscr{p}\in\P:T\in\T(\mathscr{p})\}|\ge C_1K_1 K_2\log(1/\delta)$. By Markov's inequality, we have $|\T_{\bad}|\lesssim\frac{M|\P|}{C_1K_1 K_2\log(1/\delta)}$.

Now, suppose that the lemma were false. In this case, there exists a subset $|\P'|\ge\frac{1}{2}|\P|$ and a subset $\T'(\mathscr{p})\subset\T(\mathscr{p})\cap\T_{\bad}$ for each $\mathscr{p}\in\P'$ with $|\T'(\mathscr{p})|\ge\frac{1}{2}|\T(\mathscr{p})|$. In this case, $|\T'(\mathscr{p})|$ is a $(\delta,s,pK)$-set. We now find a lower bound for $|\T'|$, where $\T'=\cup_{\mathscr{p}\in\P}\T'(\mathscr{p})$. Let 
$$J=|\{(\mathscr{p}_1,\mathscr{p}_2,T):\mathscr{p}_1,\mathscr{p}_2\in\P',T\in\T',T\in\T'(\mathscr{p}_1)\cap\T'(\mathscr{p}_2)\}|.$$
Then we have
$$J\lesssim\sum_{\mathscr{p}\in\P'}\sum_{\substack{w\in p^{-\N}\\\delta<w<\Delta}}K_2M(\delta/w)^s\cdot K_1\cdot(w/\delta)^s+|\P|M\lesssim MK_1 K_2\log(1/\delta)\cdot|\P|.$$
Here, there are $\le K_1(w/\delta)^s$ many $\delta$-balls at distance $w$ from $\mathscr{p}$, and there are $\lesssim K_2 M(\delta/w)^s$ many $\delta$-tubes in $\T'(\mathscr{p})\cap \T'(\mathscr{q})$ with $\dist(\mathscr{p},\mathscr{q})=w$. Thus, applying Cauchy-Schwarz, 
$$\frac{(M|\P|)^2}{|\T'|}\lesssim MK_1 K_2\log(1/\delta)\cdot|\P|,$$
so $|\T'|\gtrsim(K_1 K_2\log(1/\delta))^{-1}M|\P|$. Since $\T'\subset\T_{\bad}$, this contradicts our upper bound $|\T_{\bad}|\lesssim\frac{M|\P|}{C_1K_1 K_2\log(1/\delta)}$ if $C_1$ is sufficiently large.
\end{proof}

\begin{proof}[Proof of \Cref{thm: 2-s high and s low}, based on \Cref{lem: large incidence tube}]
Pick $\delta_0=\delta_0(\epsilon)$, to be chosen sufficiently small later. If $\delta>\delta_0$, then the trivial bound $|\T|_\delta\ge M$ would suffice when the implicit constant in $\gtrapprox_\epsilon$ is sufficiently large depending on $\delta_0(\epsilon)$. From now on, we assume $\delta<\delta_0$.

First, if $\delta^{8\epsilon}<\Delta$, then $\P$ is a $(\delta,s,\delta^{-17\epsilon})$ Katz-Tao set, since for $\rho\ge\Delta$ we have
$$|\P\cap B(x,\rho)|\lesssim(\rho/\Delta)^2\sup_x|\P\cap B(x,\Delta)|\lesssim\Delta^{-2}\cdot\delta^{-\epsilon^2}(\Delta/\delta)^s\le\delta^{-17\epsilon}\cdot(\rho/\delta)^s.$$
Hence, by \Cref{lem: many good tubes}, there exists a subset $\P'\subset\P$ with $|\P'|\ge\frac{1}{2}|\P|$ and for each $\mathscr{p}\in\P$ a subset $\T'(\mathscr{p})\subset\T(\mathscr{p})\in M$ such that for each $T\in\T':=\cup_{\mathscr{p}\in\P'}\T'(\mathscr{p})$,
$$|\{\mathscr{p}\in\P':T\in\T'(\mathscr{p})\}|\le C_1\delta^{-34\epsilon}\log(1/\delta).$$
Therefore,
$$|\T|\cdot\sup_{T\in\T'}|\{\mathscr{p}\in\P: T\in\T'(\mathscr{p})\}|\ge|\P'|\sup_{\mathscr{p}\in\P'}|\T'(\mathscr{p})|\gtrsim|\P|M\Longrightarrow|T|\gtrsim_\epsilon|\P|M\delta^{35\epsilon}.$$
Thus, assume $\delta^{8\epsilon}>\Delta$. In particular, $\Delta<\delta_0^{8\epsilon}$.

Next, if $\delta>\frac{\Delta}{100}$, then $\P$ is a $(\delta,2-s.\delta^{-\epsilon})$-set for analogous reasons, and we may conclude $|\T|\gtrsim_\epsilon\delta^{-1+4\epsilon}M$ by \Cref{prop: multi-grid}. Thus, we assume $\delta<\frac{\Delta}{100}$.

Now, we apply \Cref{lem: many good tubes} to the configuration $(\P,\T)$ to get a strong incidence bound. Apply \Cref{lem: pigeonhole} to $\P$, still denote as $\P$ is $\{\Delta_j\}_{j=1}^{\epsilon^{-1}}$-uniform where $\Delta_j=\Delta^{\epsilon j}$. Then we apply \Cref{prop: new scale Delta} to $(\P,\T)$ with parameter $\Delta$ to obtain a refinement $(\P_1,\T_1)$ and the corresponding $(\Delta,s,C_\Delta,M_\Delta)$-nice configuration $(\P_1^\Delta,\T_1^\Delta)$ with $C_\Delta\approx_\Delta\delta^{-\epsilon^2}$. By \Cref{prop: multi-grid} with $n=4/\epsilon$, and $C_0:=\delta^{-\epsilon^2}<\Delta^{-\epsilon/8}$, there exists a refinement $(\mathcal{Q}^\Delta,\T^\Delta)$ of $(\P_1^\Delta,\T_1^\Delta)$ such that $|\T^w\cap\D_w(\P_1)|\lessapprox_{\Delta,\epsilon}\Delta^{-3\epsilon/4}\cdot w|\D_w(\P_1)|$ for all $w\in\{1,\Delta^\epsilon,\ldots,\Delta\}$ and $\mathbf{T}\in\T^\Delta$. Now we define $\P_2=\cup_{\mathbf{p}\in\mathcal{Q}_\Delta}$ and for each $\mathscr{p}\in\P_2$, let $\T_2(\mathscr{p})=\{T\in\T_1(\mathscr{p}): T^\Delta\in\T^\Delta(\mathbf{p})\}$ and we get from \Cref{prop: new scale Delta} that
$$\sum_{\mathscr{p}\in\P_2}|\T_2(\mathscr{p})|\approx_\delta|\P|_\delta\cdot M.$$
By pigeonholing, we can find a configuration refinement $(\P_3,\T_3)$ of $(\P_2,\T_2)$ that for all $\mathscr{p}\subset\P_3$, $$\D_\Delta(\T_3(\mathscr{p})\subset\T^\Delta(p^\Delta)).$$
By \Cref{lem: many good tubes} on each $Q\in\D_\Delta(\P_3)$ and assuming $\epsilon$ is small enough, we can find a refinement $(\P_4,\T_4)$ of $(\P_3,\T_3)$ such that for all $T\in\T_4$ and each $Q\in\D_\Delta(\P_4)$,
$$|\{\mathscr{p}:Q\cap\P_4:T\in\T_4(\mathscr{p})\}|\lesssim\delta^{-\epsilon}.$$
By refining $\P_4$ further, we may assume $\P_4$ is $\{\Delta_j\}_{j=1}^{\epsilon^{-1}}$-uniform. For each $T\in\T_4$, there exists $a(T)\in[1,\delta^{-1}]\cap p^{\N}$ such that
$$|\{Q:\in\D_\Delta(\P_4):\exists Q\cap\P_4,T\in\T_4(\mathscr{p})\}|\sim a(T).$$
Also, there exists $b(T)\in[1,\delta^{-1}]$ such that there are $a'(T)\approx_\delta a(T)$ many $Q\in\D_\Delta(\P_4)$ with
$$|T\cap\P_4\cap Q|\sim b(T).$$
By pigeonholing again, there exists a pair $(a,b)$ with $a,b\in p^{\N}$ such that
$$\T_4':=\{T\in\T_4: a(T)\sim a, b(T)\sim b\}$$
preserves the number of incidences
$$|\{(p,T)\in\P_4\times\T_4':T\in\T_4(\mathscr{p})\}|\gtrapprox_\delta|\P|_\delta M.$$
Following the same estimate of $|\T_4'|$ from the proof of \cite[Proposition 4.6]{ren2023furstenberg}, we end the proof by the estimate
$|\T_4'|\gtrapprox_{\delta,\epsilon}\min\{M|\P|_\delta,M^{3/2}|\P|_\delta^{1/2},\delta^{-1}M\}.$
\end{proof}

\subsection{Proof of \Cref{lem: large incidence tube}}

The goal of this subsection is to prove the following theorem, which restates \Cref{lem: large incidence tube} and may have independent interest.

\begin{thm}\label{thm: upper bound of rich}
Let $\delta\le\Delta\in p^{-\N},\epsilon\in\frac{1}{\N}$, and let $\P$ be a multi-set of $\delta$-tubes such that for all $1\le k\le \epsilon^{-1}$, each $Q\in\D_{\Delta^{k\epsilon}}(\P)$ contains about the same number of cubes in $\P$ (with multiplicity). For $a\ge 2$ and $ab>\delta^{1-2\epsilon}|\P|$, let $\T_{a,b}$ be a set of distinct $\delta$-tubes satisfying
\begin{enumerate}
\item \label{item: good} $|T^\rho\cap\D_\rho(\P)|\le\Delta^{-\epsilon}\cdot \rho|\D_\rho(\P)|$ for all $\rho\in\{\Delta^\epsilon,\Delta^{2\epsilon},\ldots,\Delta\}$;
\item \label{item: a}$N_{\Delta,b}(T)\ge a$.
\end{enumerate}
Then $|\T_{a,b}|\lesssim_\epsilon\frac{|\P|^2}{a^3b^2}\delta^{-5\epsilon}$.
\end{thm}

The following lemma is crucial for our proof, which provides a less stronger estimate for the weaker setting that $\P$ is a $(\Delta,1,K)$-set.

\begin{lemma}\label{lem: weaker rich line}
Let $\delta\le\Delta\in p^{-\N}$. Suppose $\P\subset\P_\delta$ is a $(\Delta,1,K)$-set of (not-necessarily distinct) $\delta$-tubes, i.e., for all $\Delta\le r\le 1$, we have $|\P\cap Q|\le Kr|\P|$ for all $Q\in\D_r(\P)$. Then 
$$\sum_{\substack{T\in\T_\delta:\\N_{\Delta,b}(T)\ge 2}}N_{\Delta,b}(T)^2\lesssim K\log(1/\Delta)\cdot\frac{|\P|^2}{b^2}.$$
Consequently, if $T\subset\T_\delta$ is a set of distinct dyadic $\delta$-tubes with $N_{\Delta,b}(T)\ge a$ where $a\ge 2$, then
$$|\T_{a,b}|\lesssim K\log(1/\Delta)\cdot\frac{|\P|^2}{a^2b^2}.$$
\end{lemma}

\begin{proof}
The proof technique follows similarly from \Cref{lem: many good tubes}. Denote
$$J=|(\mathscr{p}_1,\mathscr{p}_2,\T)\in \P^2\times\T_\delta: \mathscr{p}_1^\Delta\ne\mathscr{p}_2^\Delta, \mathscr{p}_1\in T,\mathscr{p}_2\in T|.$$
On the one hand, for each $T$ with $N_{\Delta,b}(T)=a\ge 2$, it contributes $\ge ab\cdot(a-1)b\ge\frac{1}{2}a^2b^2$ many triples to $J$, so 
\begin{equation}\label{eq: contribution of each tube}
J\ge\frac{1}{2}b^2\sum_{T\in\T_\delta}|N_{\Delta,b}(T)|^2.
\end{equation}
On the other hand, for any scale $\rho\in[\Delta,1]$ and $\mathscr{p}_1,\mathscr{p}_2$ with distance $\rho$, there are $\le\rho$ many distinct $\delta$-tubes $T\in\T_\delta$ that cover both $\mathscr{p}_1$ and $\mathscr{p}_2$. Therefore, we have
\begin{equation}\label{eq: contribution of each cube pair}
J\lesssim\sum_{\mathscr{p}_1\in\P}\sum_{\rho\in[\Delta,1]\cap p^{-\N}}K|\P|\cdot\rho\cdot\frac{1}{\rho}\sim K\log(1/\Delta)\cdot|\P|^2.
\end{equation}
The estimates \eqref{eq: contribution of each tube} and \eqref{eq: contribution of each cube pair} together give the proof.
\end{proof}

\begin{proof}[Proof of \Cref{thm: upper bound of rich}]
\noindent{\bf Step 1. A weak estimate.}
Before anything else, let us check that $\P$ satisfies the condition of \Cref{lem: weaker rich line}. For any $T\in \T_{a,b}$ and $\rho=\Delta^{k \eps}$, $1\le k\le \eps^{-1}$,  we have 
$$1\leq N_{\rho,b}(T)\leq \Delta^{-\eps} \cdot \rho|\D_\rho(\P)|.$$ 
Therefore, $|\D_\rho(\P)|\geq \Delta^{\eps} \cdot \rho^{-1}$. Since each $Q\in \D_{\Delta^{k\eps}}(\P)$ contains about the same number of balls in $\P$, we obtain
\begin{equation*}
|\P \cap Q| \sim  \frac{|\P|}{|\D_\rho (\P)|} \le \Delta^{-\eps} \cdot \rho |\P|.
\end{equation*}
The same holds for general $r\in[\Delta,1]\cap p^{-\N}$ at a cost of a $\Delta^\epsilon$ factor (round $r$ up to the nearest $\Delta^{k\epsilon}$). Thus, the assumptions of \Cref{lem: weaker rich line} are satisfied by taking $K=\Delta^{2\epsilon}$, so 
\begin{equation}\label{eq: weak bound}
|\T_{a,b}|\lesssim_\epsilon\Delta^{-3\epsilon}\cdot\frac{|\P|^2}{a^2b^2}.
\end{equation}
This bound is not strong enough in general, but it  enables us to dispose of the special edge case: Pick $\Delta_0 = \Delta_0 (\eps) > 0$ a small constant that we will choose later. If $\Delta\ge\Delta_0$, then to ensure $\T_{a,b}\ne\emptyset$, we must require $a \le \Delta^{-1} \le \Delta_0^{-1}$. Thus, the estimate \eqref{eq: weak bound} will suffice upon absorbing $\Delta_0^{-1}$ into the implicit constant in $\lesssim_{\eps}$.

\noindent{\bf Step 2. Base cases.} For a fixed $\Delta<\Delta_0$, we induct on $\delta=\Delta p^{-n}$. There are two simple base cases. The first case is when $\delta>\Delta^{1+\epsilon}$ is not much smaller than $\Delta$. We will prove $|\T_{a,b}|=0$, i.e., the requirements \eqref{item: good} and \eqref{item: a} cannot both be satisfied. Indeed, if $T\in\T_{a,b}$, then $|\P\cup Q|\ge b$ for some $Q\in\D_\Delta(\P)$. By uniformity of $\P$ at scale $\Delta$, we get $\P\gtrsim b|\D_\Delta(\P)|$. Therefore, applying \eqref{item: a}, we have
$$N_{\Delta,b}(T)\ge a>\delta^{1-2\epsilon}b^{-1}|\P|\gtrsim\delta^{1-2\epsilon}|\D_\Delta(\P)|.$$
On the other hand, by \eqref{item: good}, we have
$$N_{\Delta,b}(T)\le |T^\Delta\cap \D_\Delta(\P)|\le\Delta^{-\epsilon}\cdot\Delta|\D_\Delta(\P)|=\Delta^{1-\epsilon}|\D_\Delta(\P)|.$$
However, $\delta^{1-2\epsilon}>\Delta^{1-\epsilon-2\epsilon^2}>\Delta^{1-\epsilon}$, which makes a contradiction. Thus, we have $|\T_{a,b}|=0$.

The second case is when $a$ is small. If $a\le 2\Delta^{-2\epsilon}$, then \eqref{eq: weak bound} is already a desired bound.

\noindent{\bf Step 3. Inductive step and choosing an intermediate scale.} From now on, we can assume $\delta<\Delta^{1+\epsilon},a>2\Delta^{-2\epsilon}$, and $\T_{a,b}\ne\emptyset$. For the inductive step, assume that this is already true for $\delta>\Delta p^{-n+1}$, and consider $\delta=\Delta p^{-n}$. We first choose an intermediate scale $w$ for later use of the high-low method. Let $k$ be the largest integer $1\le k\le \epsilon^{-1}$ such that
\begin{equation}\label{eq: critical scale}
a>2\Delta^{-\epsilon}\cdot w|\D_w(\P)|
\end{equation}
for $w=\Delta^{k\epsilon}$. This $k$ must exist because when $k=1$, this corresponds to our assumption $a>2\Delta^{-2\epsilon}$. Thus $k\ge 1$. We also have $k<\epsilon^{-1}$, since otherwise $a>2\Delta^{1-\epsilon}\cdot|\D_\Delta(\P)|$, contradicting \eqref{item: good}.

Since $k$ is maximal, we have the following properties.

\begin{enumerate}
\item\label{item: conc1} $a<2\Delta^{-\epsilon}\cdot(w\Delta^\epsilon)|\D_{w\Delta^\epsilon}(\P)|\le 2\Delta^{-2\epsilon}\cdot w|\D_w(\P)|$, since each $Q\in\D_w(\P)$ contains $\lesssim\Delta^{-2\epsilon}$ elements of $\D_{w\Delta^{\epsilon}}(\P)$. Therefore, we obtain the bound
\begin{equation}\label{eq: bound of a}
2\Delta^{-\epsilon}\cdot w|\D_w(\P)|<a<2\Delta^{-2\epsilon}\cdot w|\D_w(\P)|.
\end{equation}
\item\label{item: conc2} For any $Q\in\D_w(\P)$, we claim that $\P^Q=\phi_Q(\P\cap Q)$ is a $(\Delta/w,1,\Delta^{-\epsilon})$-set, where $\phi_Q$ is the affine transform that maps $Q$ to $\Z_p^2$. This is because for $\rho=\Delta^{m\epsilon},m>k$, we have $w|\D_w(\P)|\le\frac{a}{p\Delta^{-\epsilon}}\le\rho|\D_\rho(\P)|$, and every $\rho$-tube in $\D_\rho(\P)$ contains about the same number of cubes in $\P$. Thus, for any $\tilde Q\in\D_\rho(\P)$ and $Q\in\D_w(\P)$,
$$|\tilde Q\cap\P|\sim\frac{|\P|}{|\D_\rho(\P)|}\le\frac{|\P|}{|\D_w(\P)|}\cdot\frac{\rho}{w}\sim|Q\cap\P|\cdot\frac{\rho}{w}.$$
By rounding $r$ to the nearest $\Delta^{m\epsilon}$, we can ensure that 
$$|\tilde Q\cap P|\le\Delta^{-\epsilon}|Q\cap\P|\cdot\frac{\rho}{w}$$
for all $\tilde Q\in\D_r(\P),\Delta\le r\le w$. 
\end{enumerate}
For the rest of the proof, we will only use the above properties \eqref{item: conc1} and \eqref{item: conc2} of the scale $w=\Delta^{k\epsilon}$.

\noindent{\bf Step 4. Incidences and high-low method.}
For each $Q\in\D_w(\P)$, let $\U_Q=\{T\cap Q: T\in\T_{a,b}\}$ be a set of $\delta\times w$-tubelets, and $\U=\bigsqcup_{Q\in\D_w(\P)}$. For $u\in\U$, define $N_{\Delta,b}(u)$ as the set number of $Q\in\D_\Delta(\P)$ such that $|u\cap Q\cap\P|\ge b$, and $m(u)$ as the number of $\delta$-tubes $T\in\T_{a,b}$ that contains $u$. We assign $u$ with the weight $N_{\Delta,b}(u)$ and consider weight incidences. With this convention, we have
$$a|\T_{a,b}|\le I(\U,\T_{a,b}):=\sum_{u\in\U}N_{\Delta,b}(u)\cdot|\{T\in\T_{a,b}:u\subset T\}|\ge\sum_{u\in\U}N_{\Delta,b}(u)m(u).$$
On the other hand, by \eqref{eq: critical scale} and \eqref{item: good}, 
$$\frac{a}{2}|\T_{a,b}|\ge\sum_{T\in\T_{a,b}}|T^w\cap\D_w(\P)|\cdot|\T_{a,b}|\ge\sum_{u\in\U}m(u).$$
Therefore, at least half of the incidences between tubes and tubelets involve tubelets with $N_{\Delta,b}(u)\ge2$. Now, let us replace $\U$ by the tubelets $u\in\U$ with $N_{\Delta,b}(u)\ge 2$. 

For each $\mathbf{T}\in\T_{\delta/w}$, denote by $\U_\mathbf{T}=\{T\cap Q: T\in\T_{a,b},T\subset\mathbf{T}\}$ the set of tubelets that are intersections of $Q\in\D_w(\P)$ and $T\in\T_{a,b}$ that are contained in $\mathbf{T}$. After rescaling $\mathbf{T}$ to $\Z_p^2$, the tubelets in $\U_\mathbf{T}$ become $w$-cubes and the $\delta$-tubes in $\T_\mathbf{T}:=\T_{a,b}\cap \mathbf{T}$ become $w$-tubes. This allows us to decompose the weighted incidence count as 
$$\I(\U,\T_{a,b})=\sum_{\mathbf{T}\in\T_{\delta/w}}\I(\U_\mathbf{T},\T_\mathbf{T}).$$
Applying the high-low estimate in \Cref{prop: High-low method} to a thickening index $S = \Delta^{-\eps/100} \in (w^{-\eps/100}, w^{-1})$, we have
\begin{align}
\begin{split}
a|\T_{a,b}|&\le \I(\U,\T_{a,b})\\
&=\sum_{\mathbf{T}\in\T_{\delta/w}}\I(\U_\mathbf{T},\T_\mathbf{T})\\
&\le S^{1/2} w^{-1/2}\sum_{\mathbf{T}\in\T_{\delta/w}}|\T_\mathbf{T}|^{1/2}\left(\sum_{u\in\U_\mathbf{T}}N_{\Delta,b}(u)\right)^{1/2}+S^{-1}\sum_{\mathbf{T}\in\T_{\delta/w}}\I(\U_\mathbf{T},\T_\mathbf{T}^{S\delta}),
\end{split}
\end{align}
where $\T_\mathbf{T}^{S\delta}:=\{T^{S\delta}:T\in\T_\mathbf{T}\}$ is counted with multiplicity. We denote by
$$(A):=\sum_{\mathbf{T}\in\T_{\delta/w}}|\T_\mathbf{T}|^{1/2}\left(\sum_{u\in\U_\mathbf{T}}N_{\Delta,b}(u)\right)^{1/2}$$
the high-frequency term, and
$$(B):=\sum_{\mathbf{T}\in\T_{\delta/w}}\I(\U_\mathbf{T},\T_\mathbf{T}^{S\delta})$$
the low-frequency term.

\noindent{\bf Step 5. High-frequency case.}
Using Cauchy-Schwarz and the fact that $N_{\Delta,b}\ge 2$ for any $u\in\U$, the high-frequency term $(A)$ is bounded by
\begin{align}
\begin{split}
(A)&\le\left(\sum_{\mathbf{T}\in\T_{\delta/w}}|\T_\mathbf{T}|\right)^{1/2}\left(\sum_{\mathbf{T}\in\T_{\delta/w}}\sum_{u\in\U_\mathbf{T}}N_{\Delta,b}(u)\right)^{1/2}\\
&=|\T_{a,b}|^{1/2}\left(\sum_{Q\in\D_w(\P)}\sum_{u\in\U_Q}N_{\Delta,b}(u)\right)^{1/2}\\
&\lesssim |\T_{a,b}|^{1/2}\left(|\D_w(\P)|\cdot\Delta^{-2\epsilon}(\frac{|\P|}{b|\D_w(\P)|})^2\right)^{1/2}.
\end{split}
\end{align}
Here, the last line comes from \eqref{item: conc2} and the weak estimate in \Cref{lem: weaker rich line}. If the high-frequency term dominates, i.e., $a|\T_{a,b}|\lesssim_\epsilon w^{-1/2}(A)$, we get by setting $S=\Delta^{-\eps/100}$ and \eqref{eq: bound of a},
$$|\T_{a,b}|\lesssim\Delta^{-2\epsilon} S\frac{w^{-1}}{a^2b^2}|\D_w(\P)|\cdot(\frac{|\P|}{|\D_w(\P)|})^2\lesssim \Delta^{-5\epsilon}\frac{|\P|^2}{a^3b^2}.$$

\noindent{\bf Step 6. Low-frequency case.} Now suppose the low-frequency case dominates, i.e.,
\begin{align}
\begin{split}
Sa|\T_{a,b}|&\lesssim \sum_{\mathbf{T}\in\T_{\delta/w}}\I(\U_\mathbf{T},\T_\mathbf{T}^{S\delta})\\
&=\sum_{\mathbf{T}\in\T_{\delta/w}}\sum_{T\in\T_\mathbf{T}}\sum_{u\in\U_\mathbf{T}\cap T^{S\delta}}N_{\Delta,b}(u)\\
&=\sum_{T\in\T_{a,b}}f(T),
\end{split}
\end{align}
where $f(T):=\sum_{u\in\U_\mathbf{T}\cap T^{S\delta}}N_{\Delta,b}(u)$ for $T\in\T_\mathbf{T}$. Thus, there exists a subset $\T_{a,b}'\subset\T_{a,b}$ with $|\T_{a,b}'|\ge\frac{1}{2}|\T_{a,b}|$ such that for each $T\in\T_{a,b}'$, we have $f(T)\gtrsim Sa$.

Now, for $Q\in\D_\Delta(\P)$, we denote by $g_T(Q)$ be the number of $u\in\U_\mathbf{T}\cap T^{S\delta}$ such that $|u\cap Q\cap\P|\ge b$. Then we immediately have
$$|T^{S\delta}\cap Q\cap \P^{S\delta}|\ge g_T(Q)b,$$
where $\P^{S\delta}$ is counted with multiplicity. By double counting and the definition of $N_{\Delta,b}(u)$, we have for $T\in\T_{a,b}'$,
$$\sum_{Q\in\D_\delta(\P)}g_T(Q)=f(T)\gtrsim a.$$
We also have $g_T(Q)\le S$ since at most this many tubelets in $\U_\mathbf{T}$ can intersect $T^{S\delta}\cap Q$, which is a $S\delta\times\Delta$ rectangle. Hence, by pigeonholing, we can find a  $k_T\in[1,S]\cap p^{\N}$ and a set $\mathcal{Q}_T\subset\D_\Delta(\P)$ with $|\mathcal{Q}_T|\gtrsim \frac{Sa}{k_T \log S}$ such that $g_T(Q)\sim k_T$ for all $Q\in\mathcal{Q}_T$. By another pigeonholing, 
we can further choose a $k\in p^{\N}$ and a subset $\T_{a,b}'\subset\T_{a,b}''$ with $|\T_{a,b}'|\gtrsim (\log S)^{-1} |\T_{a,b}''|$ such that $k_T = k$ is the same for all $T\in\T_{a,b}'$. 

Consider the set $\tilde \T:=\{T^{S\delta}: T\in\T_{a,b}'\}$, where we delete the duplicate elements. It has cardinality $\gtrsim |\T_{a,b}'|$ because at most $\lesssim 1$ many $\delta$-tubes share a common ancestor $T^{S\delta}$.

We aim to apply the inductive hypothesis to $\tilde\T$ and $\tilde\P:=\{\mathscr{p}^{S\delta}:\mathscr{p}\in\P\}$ (with multiplicity, so $|\tilde\P|=|\P|$) over the scale $S\delta$, where $a'=\frac{Sa}{k \log S}$ and $b'=bk$.

We first check that $a', b'$ satisfy the conditions of Lemma \ref{lem: large incidence tube}. Indeed, $a' \gtrsim \frac{a}{\log S} \ge 2$ since $a \ge 2\Delta^{-2\eps}$. Furthermore, $a'b' \gtrsim S^{1-\eps} (\log S)^{-1} ab > \max\{ (S\delta)^{1-2\eps} |\P|, 2 \}$.

We next check that $\tilde\T$ satisfies the requirement \eqref{item: good}. Indeed, for all $T\in\T_{a,b}'$ and $\rho\in\{\Delta^{\epsilon},\Delta^{2\epsilon},\ldots,\Delta\}$,
$$|(T^{S\delta})^\rho\cap\D_\rho(\P)|=|T^\rho\cap\D_\rho(\P)|\le\Delta^{-\epsilon}|\Delta_\rho(\P)|.$$
Here, we have $(T^{S\delta})^\rho=T^\rho$ since $\rho\ge\Delta\ge S\delta$, using the fact that $\delta\le\Delta^{1+\epsilon}$.

We finally check requirement \ref{item: a}. For every $T^{S\delta} \in \tilde\T$, we have $N_{\Delta, kb} (T^{S\delta}) \ge a'$, the size of $\mathcal{Q}_T$.

Thus, we may apply the inductive hypothesis to obtain (after using $k \lesssim S$):
\begin{equation*}
    |\T_{a,b}| \lesssim S^2 \log S \cdot |\tilde\T| \lesssim_\eps S^2 \log S \cdot \frac{|\P|^2}{(a')^3 (b')^2} (S\delta)^{-5\eps} \le S^{-\eps} \cdot \frac{|\P|^2}{a^3 b^2} \delta^{-5\eps}.
\end{equation*}

Therefore, we close the inductive step and thus finish the proof.
\end{proof} 

\subsection{Returning to the discretized version of Furstenberg set}\label{subsec: Returning to the discretized version of Furstenberg set}

We will now use Theorem \ref{thm: 2-s high and s low} and Lemma \ref{lem: many good tubes} to supply the missing proofs of Theorems \ref{thm: discretized s<t<2-s}-\ref{thm: discretized s+t>2}.

\begin{proof}[Proof of \Cref{thm: discretized s<t<2-s}]
Notice that the bound of $|\T|_\delta$ in \Cref{thm: 2-s high and s low} does not rely on our choice of $\Delta$. Therefore, we can take $\Delta=(|\P|_\delta\delta^{s})^{-1/(2s-2)}\in[\delta,1]$ in \Cref{thm: 2-s high and s low}. In this way, the two bounds \eqref{eq: 2-s high} and \eqref{eq: s low} together implies $|\P\cap Q|_\delta\le\delta^{-\epsilon^2}\cdot\max\{\rho^{2-s}|\P|_\delta,(\rho/\delta)^s\}$, thus satisfy the requirement \eqref{eq: strong space condition}. 
\end{proof}

\begin{proof}[Proof of \Cref{thm: discretized 0<t<s<1}]
Note that since $\P$ is a $(\delta, t, \delta^{-\eta})$-set and $s \ge t$, it satisfies the spacing condition
$$
|\P \cap B(x, w)| \le \delta^{-\eta} |\P| w^t \le (\delta^{t-\eta} |\P|) \left( \frac{w}{\delta} \right)^t \le (\delta^{t-\eta} |\P|) \left( \frac{w}{\delta} \right)^s.
$$
Apply \Cref{lem: many good tubes} with $K_1 = |\P| \delta^{t-\eta}$, $K_2 = \delta^{-\eta}$, and $\Delta=1$ to obtain $\P'$ and $\T'$. By a double counting argument, we get
$$
|\P| M \lesssim \sum_{\mathscr{p} \in \P'} |\T'(p)| = \sum_{T \in \T'} |\mathscr{p}\in\P':T\in\T'(\mathscr{p})|\le C_1|\P| \delta^{t-2\eta} \log(1/\delta) |\T|_\delta.
$$
Rearranging gives $|\T|_\delta \gtrsim \delta^{-t+2\eta} M$, as desired.
\end{proof}

\begin{proof}[Proof of \Cref{thm: discretized s+t>2}]
Since $t \ge 2-s$, we see that $\P$ is a $(\delta,2-s,\delta^{-\epsilon^2})$-set. Set $\Delta=\delta$ in \Cref{thm: 2-s high and s low}. In this case, the bound \eqref{eq: 2-s high} is exactly our condition that $\P$ is a $(\delta,2-s,\delta^{-\epsilon^2})$-set, which implies $|\P|_\delta\gtrsim\delta^{s-2+\epsilon^2}$, and 
$$|\T|_\delta\gtrapprox_{\delta,\epsilon}\delta^{35\epsilon}\cdot\min\{M|\P|_\delta,M^{3/2}|\P|_\delta^{1/2},\delta^{-1}M\}\gtrsim_\epsilon\delta^{-1+\epsilon}M.$$
\end{proof}

\bibliographystyle{plain}
\bibliography{references.bib}

\end{document}